\documentclass[11pt,a4paper,leqno]{amsart}
\usepackage{latexsym}
\usepackage{floatflt}
\usepackage[active]{srcltx}
\usepackage{hyperref}
\usepackage{amscd}

\usepackage{amsmath, amsthm, amssymb}
\usepackage[all,cmtip]{xy}
\usepackage[parfill]{parskip}    
\usepackage{floatflt}
\usepackage[noadjust]{cite}
\def\ts{\thinspace}%
\def\eps{\varepsilon}%
\def\phi{\varphi}

\def\"#1{{\accent"7F #1\penalty10000\hskip 0pt plus 0pt}} 
\def\ts{\thinspace}%
\def\ge{\geqslant}%
\def\le{\leqslant}%

\def\CP{{\mathbb C} {\mathbb P}}

\newtheorem{thm}{Theorem}[subsection]

\newtheorem{Example}[thm]{Example}

\newtheorem{Counterexample}[thm]{Counterexample}

\newtheorem{remark}[thm]{Remark}

\newtheorem{Fact}[thm]{Fact}
\newtheorem{Nothing}[thm]{$\!\!\!$}

\newcommand{\be}{\begin{equation} }

\newcommand{\ene}{\end{equation} }

\newcommand{\ba}{\begin{eqnarray}}
\newcommand{\ea}{\end{eqnarray}}
\newcommand{\ban}{\begin{eqnarray*}}
\newcommand{\ean}{\end{eqnarray*}}

\DeclareMathOperator{\A}{A} 
\newcommand{\rightquotient}[2]{\raisebox{0.6 ex}{\ensuremath{#1}} \hspace{-0.3ex} \diagup\hspace{-0.4ex} \raisebox{-0.7 ex}{\ensuremath{#2}}}
\newcommand{\lie}[1]{\mathfrak{#1}}
\newcommand{\hor}[1]{#1^{\text{hor}}}
\newcommand{\ver}[1]{#1^{\text{vert}}}
\newcommand{\isom}{\cong}

\begin{document}

\abovedisplayskip=6pt plus3pt minus3pt \belowdisplayskip=6pt
plus3pt minus3pt
\title[Almost nonnegative curvature operator and principal bundles]
{Manifolds with Almost nonnegative curvature operator and principal bundles}

\abstract{We study manifolds with almost nonnegative curvature operator (ANCO) and
provide first examples of closed simply connected ANCO mannifolds 
that do not admit nonnegative curvature operator.}
\endabstract

\author{Martin Herrmann, Dennis Sebastian, and Wilderich Tuschmann}
\address{Fakult\"at f\"ur Mathematik
\\ Karlsruher Institut f\"ur Technologie (KIT)
\\ Kaiserstra\ss{}e 89-93
\\ D-76133 Karlsruhe, Germany}
\email{martin.herrmann@kit.edu, tuschmann@kit.edu}

\maketitle

\section{Introduction}
The present note sets out to study and construct closed smooth manifolds
with almost nonnegative curvature operator, here in short termed ``ANCO'' manifolds.
Its motivation mainly stems from the following sources.

Whereas recent work of Ch.~B\"ohm and B.~Wilking ([BW08]),
showing that manifolds with positive curvature operator are space forms, 
has, when combined with previous works on the subject, also entailed a final 
classification of manifolds with nonnegative curvature operator,
the case of almost nonnegative curvature operator is still 
far from being understood.\\
Moreover, manifolds with almost nonnegative curvature operator
equally deserve special attention.
Very recently, J.~Lott has developed a first structure theory for manifolds
collapsig under a lower bound on the curvature operator and 
extended
the results of Cheeger and Gromov on F-structures in collapse 
under both-sided bonds on sectional curvature to this setting (see [Lo12]).
His work also implies that in obtaining a more general picture
capturing all collapsed directions as in the work of Cheeger-Fukaya-Gromov on 
N-structures, ANCO manifolds will here play the same crucial role as
almost flat manifolds do there and manifolds with almost nonnegative sectional curvature
do play in collapsing with curvature bounded from below.

Our first result shows that even in the simply connected case 
one can distinguish the class of manifolds with almost nonnegative operator
from the class with nonnegative one,
and that there are in fact plenty of spaces which do so.

\noindent {\bf Theorem 1} \thinspace {\it In dimension $n=7$ and each dimension $n\ge 9$ 
there exist infinite sequences of closed simply connected $n$-manifolds of pairwise
distinct homotopy type which all admit almost nonnegative curvature operator
but none of which supports a metric with nonnegative curvature operator.}

Theorem $1$ provides, at least to our knowledge, the first examples of closed simply
connected ANCO manifolds which do not support any metric with nonnegative curvature operator.

In this regard, Theorem 1 is also in strong contrast to the fact that
until now one does not know of a single example of a closed simply connected manifold 
with almost nonnegative sectional curvature for which it could be shown that it
does not also admit  (strictly)  nonnegative sectional curvature.

The simplest examples of simply connected ANCO manifolds not admitting nonnegative
curvature operator constructed here are given
by (simply connected) total spaces of
principal circle bundles over $\CP^1 \times \CP^2$ among which there is also the non-trivial
$S^5$ bundle over $S^2$.
Interestingly enough, these
already give rise to a phenomenon which is impossible in the case of manifolds
with nonnegative curvature operator:

\noindent {\bf Theorem 2} \thinspace {\it There exist homeomorphic closed simply connected 
$7$-manifolds with almost nonnegative curvature operator which are not diffeomorphic.}

To obtain the spaces which yield the above theorems,
we first construct ANCO metrics on the total spaces of certain principal bundles
over (A)NCO manifolds (see Theorem~3.1). This is done in section $3$ using a criterion of the first-named author ([Her12]). The proofs of Theorem~$1$ and Theorem~$2$, also given there, then
follow easily
from the classification of manifolds with nonnegative curvature operator and work of
Kreck and Stolz ([KS88]).

The other parts of this work are structured as follows:
In the next section, we provide relevant definitions and preliminaries,
together with some first examples of ANCO manifolds.
We also take the opportunity to draw attention to
a Betti number estimate for ANCO manifolds
due to P.~B\'erard ([Be88]).
Section $4$ concludes this work with a number of further remarks and open questions.

The last-named author wishes to thank John Lott for communicating his work to him, 
and besides John also Anton Petrunin and Burkhard Wilking for several stimulating conversations on ANCO spaces.

\section{Preliminaries on ANCO manifolds}

\noindent {\bf Definition 2.1} \thinspace {\it A closed smooth
manifold $M$ is said to admit {\rm almost nonnegative curvature operator (ANCO) }
if for all $\eps>0$ there exists a Riemannian metric $g_\eps$ on $M$ such that
all eigenvalues $\lambda_\eps$ of the associated curvature operator 
on bivectors and the
diameter of $g_\eps$
satisfy the inequality
$$\lambda_\eps \ts \cdot \ts   diam(M,g_\eps)  \ts  ^2 \ > \ -\eps \ .$$

If, in addition, 
when replacing $-\eps$ by~$\eps$, the reverse inequality is also satisfied,
$M$ is said to have {\rm almost flat curvature operator}. 
}

\

Notice that a manifold has almost flat curvature operator if and only
if it is almost flat, i.e., diffeomorphic to an infranil manifold. 
Moreover, as in the case of nonnegative curvature operator, 
Riemannian products of ANCO manifolds are ANCO.

If one is not especially concerned about the simply connected case,
there is, as is evidenced by the following fact,
indeed a multitude of ANCO manifolds.

\noindent {\bf Fact 2.1} \thinspace 
Any finitely generated almost nilpotent group $\Gamma$
can be realized as the fundamental group of an ANCO manifold.

This can be seen, for instance, by using that products of ANCO manifolds are ANCO 
and mimicking a construction of Wilking
(see [Wil00]) to make any such $\Gamma$ act cocompactly,
freely and properly discontinuously by isometries 
on a product of a connected and simply connected nilpotent Lie group, equipped
with a left invariant metric, and a special unitary group with its bi-invariant metric.\

Of course, ANCO manifolds always admit almost nonnegative sectional curvature,
but whether the converse is also true seems presently unclear. When combining
the work in [Wil00] with Fact 2.1, one actually sees:

\noindent {\bf Fact 2.2} \thinspace
A finitely generated group $\Gamma$ 
is isomorphic to the fundamental group of an ANCO manifold
iff  $\Gamma$ is isomorphic to the fundamental group of a closed manifold with
almost nonnegative sectional curvature
iff  $\Gamma$ is isomorphic to the fundamental group of a closed manifold with
almost nonnegative Ricci curvature.

For a brief survey on principal obstructions to admitting metrics of almost nonnegative
sectional or Ricci curvature we refer to [KPT10], and [Tu12] also describes
detailed examples of such manifolds. 
Many examples of collapsing with curvature operator bounded from below 
can be found in Lott's recent work [Lo12]. The only further obstruction to admitting
ANCO we are aware of is the following Betti number estimate which is
contained, somewhat hidden, in work of P.~B\'erard (see [Be88]).
That it also holds for manifolds of (almost) nonnegative sectional curvature
is a long-standing open conjecture by Gromov, and for NCO manifolds it just follows
from the fact that there harmonic forms are always parallel.

\noindent
{\bf Proposition 2.3 (B\'erard)} \thinspace {\it If $M^n$ admits ANCO, 
then $b_k(M^n)\le  \binom{n}{k}$.}

B\"ohm and Wilking have shown that manifolds with positive
curvature operator are space forms (see [BW08]). Together with earlier
results on the classification of manifolds with nonnegative curvature operator,
this yields a complete geometric description of all such manifolds (compare, e.g.,
[CLN06] and the further references given there). Specializing this classification
to the simply connected case (which is the only one we need here) gives

\noindent {\bf Proposition 2.4} \thinspace {\it 
A closed simply connected manifold with nonnegative curvature operator
is isometric
to a Riemannian product of

(a) standard spheres with metrics of nonnegative curvature operator;

(b) closed K\"ahler manifolds biholomorphic to complex projective spaces
with nonnegative curvature operator;

(c) compact irreducible Riemannian symmetric spaces with their natural metrics of nonnegative curvature operator.
}

\

\noindent {\bf Remark 2.5} \thinspace
In particular, as follows from Cartan's classification of symmetric spaces of compact type
(see, e.g., [Hel01]), in each fixed dimension there exist thus only a finite number of 
diffeomorphism types of closed simply connected manifolds with nonnegative curvature operator.

\section{Almost Nonnegative Curvature Operator \\ on Principal Bundles}

Here we discuss the curvature operator 
of a principal fibre bundle $\pi: P \to M$ over a compact manifold $M$ 
whose fibre is a compact Lie group~$G$.

In [FY92] K.~Fukaya and T.~Yamaguchi showed that the total spaces of such bundles admit almost nonnegative sectional curvature, if the base space $M$ does so.

To prove this result they used metrics of the following type. Let $\gamma$ be a connection form on $P$ , $b$ a biinvariant metric on $G$ and $g^M$ a Riemannian metric on $M$. Then for every $t>0$
\[g^t(X,Y)=g^M(\pi_*X,\pi_* Y)+t^2 b(\gamma(X),\gamma(Y)), \quad X,Y \in \mathrm{T}_pP,\]
defines a Riemannian metric on $P$ for which $\pi$ is a Riemannian submersion with totally geodesic fibres.

Let $\tilde{P}=\rightquotient{P}{[G,G]}$, so that there is the commutative diagram
\[\xymatrix@R=0.2cm{P\ar[dd]_{\pi} \ar[dr]^{\pi_1}&\\&\tilde{P}\ar[dl]^{\pi_2}\\M&}\]
where $\pi_1:P\to \tilde{P}$ is a principal bundle with fibre $[G,G]$ and $\pi_2: \tilde{P} \to M$ is a principal bundle with fibre the abelian Lie group
$\rightquotient{G}{[G,G]}$.

As $g^t$ is invariant under the action of $G$ on $P$, there is an induced metric $\tilde{g}^t$ on $\tilde{P}$ which makes  $\pi_1$ and $\pi_2$ also into Riemannian submersions.

The connection form $\gamma$ induces connection forms $\gamma_1$ and $\gamma_2$ for $\pi_1$ and $\pi_2$, taking values in $[\lie{g},\lie{g}]$ and $\rightquotient{\lie{g}}{[\lie{g},\lie{g}]}$, respectively.

If -- slightly abusing notation -- we denote the biinvariant metrics induced by $b$ on $[\lie{g},\lie{g}]$ and $\rightquotient{\lie{g}}{[\lie{g},\lie{g}]}\isom [\lie{g},\lie{g}]^\top$again by $b$, then the following holds:

\[g^t(X,Y)=\tilde{g}^t({\pi_1}_*X,{\pi_1}_* Y)+t^2 b(\gamma_1(X),\gamma_1(Y)), \quad X,Y \in \mathrm{T}_pP,\]
and
\[\tilde{g}^t(X,Y)=g^M({\pi_2}_*X,{\pi_2}_* Y)+t^2 b(\gamma_2(X),\gamma_2(Y)), \quad X,Y \in \mathrm{T}_p\tilde{P}.\]

The following result now gives a criterion for $P$ to admit almost nonnegative curvature operator with respect to these metrics $g^t$. 

\noindent
{\bf Theorem 3.1} \thinspace {\it
Let $(M,g^M)$ have nonnegative curvature operator. Then $(P,g^t)$ has almost nonnegative curvature operator as $t\to0$, iff $(P, g^t)$ is locally isometric to $(\tilde{P}\times [G,G], \tilde{g}^t \times t^2 b)$ for every $t$. 
}

In the case where $M$ only admits almost nonnegative curvature operator, we still have one of the implications.

{\bf Theorem 3.2} \thinspace {\it
Let $M$ admit ANCO. If $G$ is abelian, then $P$ has ANCO.
}

Notice that 
Theorems 3.1 and 3.2 
cannot be 
directly used to obtain almost nonnegative curvature operator on associated non-principal bundles, since lower bounds for the 
the eigenvalues of the curvature operator are, in general, not preserved
 under Riemannian submersions.

Both the above theorems
follow from the next proposition, contained in [Her12], and Theorem 3.2 independently
also follows  from the calculations in [Se11] and [Lo12].
At the end of this section,
we will use Theorem 3.2 to prove Theorems 1 and 2.

{\bf Main Proposition 3.3} \thinspace {\it
Let $\Omega$ denote the curvature form of $\gamma$.
\begin{enumerate}
\item Let $(M,g^M)$ have nonnegative curvature operator. Then $(P,g^t)$ has almost nonnegative curvature operator as $t \to0$ iff  \mbox{$\mathrm{Im}(\Omega) \subset [\mathfrak{g},\mathfrak{g}]^\top$.}
\item Let $M$ admit almost nonnegative curvature operator. If $\mathrm{Im}(\Omega) \subset [\mathfrak{g},\mathfrak{g}]^\top$, then $P$ admits almost nonnegative curvature operator.  \end{enumerate}
}
\vspace{-0.2cm}
\begin{proof}[Proof of theorems 3.1 and 3.2 using proposition 3.3.]
Assume that $(M,g^M)$ has nonnegative curvature operator. By the proposition, $(P,g^t)$ has almost nonnegative curvature operator as $t \to 0$ if and only if the image of the curvature form of $\gamma$ is contained in the orthogonal complement of $[\mathfrak{g},\mathfrak{g}]$. This means that the curvature form $\Omega_1$ of $\gamma_1$ is zero. Therefore $\gamma_1$ is locally isomorphic to the canonical connection form on $\tilde{P}\times [G,G]$. By the definition of the metrics $g^t$ and $\tilde{g}^t$  this gives locally defined isometries.

The other implication is given by the fact, that $\tilde{P}$ satisfies the condition of Proposition 3.3(1). Therefore $(\tilde{P},\tilde{g}^t)$ has almost nonnegative curvature operator as $t\to 0$ and so the local Riemannian product $P$. Theorem 3.2 is a direct consequence of 3.3~(2), since in this case $[\mathfrak{g},\mathfrak{g}]=\{0\}$.
\end{proof}\vspace{-0.2cm}
The proof of the above proposition is given at the end of this section
and is based upon the following preparations.
Let $\nabla^t$ denote the Levi-Civita connection of $g^t$ and $\nabla:=\nabla^1$. Furthermore, let $\mathrm{A^t}$ denote the O'Neill-tensor of $\pi$ with respect to $g^t$ and  $A:=A^1$.
The Koszul formula implies

{\bf Lemma 3.4 \thinspace}(see [FY92], Lemma 2.3) \thinspace {\it
For horizontal vector fields $X,Y$ and vertical vector fields $V,W$ we have:
\begin{enumerate}
\item $\nabla^t_X Y=\nabla_X Y,\qquad \nabla^t_V W= \nabla_V W$ 
\item $\A^t_X Y=\A_X Y, \qquad\A^t_X V= t^2 \A_X V$
\item $\ver{(\nabla^t_X V)}=\ver{(\nabla_X V)}=[X,V]$
\item $\hor{(\nabla^t_V X)}=t^2 \hor{(\nabla_V X)}$
\end{enumerate}
}\vspace{-0.4cm}

{\bf Lemma 3.5} \thinspace {\it If $X$, $Y$ are horizontal vector fields and $V$, $W$ are vertical vector fields, then:
\begin{align*}g^t((\nabla^t_W \A^t)_X Y,V)&=-t^2 \tfrac{1}{4}b(\gamma([X,Y]),[V,W]) \\ &\qquad-t^4(g(\A_{\nabla_W X} Y,V)+g(\A_X(\nabla_W Y), V)).\end{align*}}\vspace{-0.6cm}

\begin{proof}
We have
\begin{align*}g^t((\nabla^t_W \A^t)_X Y,V)&=t^2 g(\nabla_W(\A_X Y), V)\\ &\qquad-t^4(g(\A_{\nabla_W X} Y,V)+g(\A_X(\nabla_W Y), V)).\end{align*}
Furthermore $\A_X Y=\tfrac{1}{2}\ver{[X,Y]}$. From the Koszul formula it follows that
 \[t^2 g(\nabla_W(\A_X Y), V)=t^2\tfrac{1}{2} g(\nabla_W(\ver{[X,Y]}), V)=-t^2 \tfrac{1}{4}b(\gamma([X,Y]),[V,W]).\]
\end{proof}

We can now use this lemma together with the O'Neill formulas to calculate the components of the curvature operator of $(P,g^t)$. To do this, let us fix a point $p \in P$ and an orthonormal basis (w.r.t. $g^1$) 
\[X_1, \dots, X_n, X_{n+1},\dots, X_{n+r}\]
of $\mathrm{T}_pM$, where $n=\dim M$, $r=\dim G$, $X_i$, $i\leq n$, is horizontal and $X_a$, $a\geq n+1$, is vertical. We extend $X_1, \dots , X_n$ to vector fields which locally form a basis of the horizontal space and use the same symbols to denote these vector fields. 
For $a \geq n+1$ we extend $X_a$ as the canonical vector field given by the right action of $G$, so $X_a(q)=\tfrac{\mathrm{d}}{\mathrm{d}t}\big\vert_{t=0}\left(q.\exp(t \gamma(X_a))\right)$.


{\bf Proposition 3.6} \thinspace {\it
For the curvature tensor $R^t$ of $P$ with respect to the metric $g^t$ we have
\begin{align*}
R^t(X_i,X_j,X_l,X_k)&=R_M(X_i,X_j,X_l,X_k)+t^2\big(2g(\A_{X_i}X_j,\A_{X_l}X_k)\\
	&\qquad+g(\A_{X_i} X_l, \A_{X_j}X_k)-g(\A_{X_j} X_l, \A_{X_i}X_k)\big),\\
R^t(\tfrac{X_a}{t},\tfrac{X_b}{t},\tfrac{X_d}{t},\tfrac{X_c}{t})&=\tfrac{1}{t^2} R_{(G,b)}(X_a,X_b,X_d,X_c),\\
R^t(X_i,X_j,\tfrac{X_b}{t},\tfrac{X_a}{t}) &=- \tfrac{1}{2}b(\gamma([X_i,X_j]),[X_a,X_b]) \\
	&\qquad + t^2 \Big(g(\A_{\nabla_{X_a} X_i} X_j,X_b)+g(\A_{X_i}(\nabla_{X_a} X_j), X_b)\\
	&\qquad\qquad-g(\A_{\nabla_{X_b} X_i} X_j,X_a)-g(\A_{X_i}(\nabla_{X_b} X_j), X_a)\\
	&\qquad\qquad+g(\A_{X_i} X_b, \A_{X_j} X_a)\Big)\\
R^t(X_i,X_j,\tfrac{X_a}{t},X_k) &=-t\Big(g(\nabla_{X_k} (\A_{X_i} X_j) , X_a )-  g(\A_{\hor{(\nabla_{X_k} X_i)}} X_j , X_a )\\
&\qquad\qquad- g(\A_{X_i}(\hor{(\nabla_{X_k} X_j)}) , X_a )\Big),\\
R^t(\tfrac{X_a}{t},\tfrac{X_b}{t},\tfrac{X_c}{t},X_i) &=0,\\
R^t(X_i,\tfrac{X_a}{t}, \tfrac{X_b}{t}, X_j)&=-\tfrac{1}{4}b(\gamma([X_i,X_j]),[X_a,X_b])\\
&\qquad-t^2\Big(g(\A_{\nabla_{X_b} X_i} X_j,X_a)+g(\A_{X_i}(\nabla_{X_b} X_j), X_a)\\
&\qquad\qquad+g(\A_{X_i} X_a,\A_{X_j}X_b)\Big)
\end{align*}
where $i,j,k,l\leq n$ and $n<a,b,c,d \leq n+r$.
}

\begin{proof} This follows from the above lemmata by a straightforward calculation.

\end{proof}

Hence, in the basis \[X_i \wedge X_j, i<j\leq n, \;\tfrac{1}{t^2} X_a \wedge X_b, n <a<b, \;\tfrac{1}{t}X_i \wedge X_a, i\leq n<a,\]
the curvature tensor is given by a matrix
\[\left(\renewcommand{\arraystretch}{1.4}
\begin{array}{c|c|c}\hat{R}_M  & B & 0 \\\hline  B^\top& \tfrac{1}{t^2} \hat{R}_{(G,b)}&0\\\hline 0 & 0 & A \end{array}\renewcommand{\arraystretch}{1}\right)+ t \;C_1+t^2 \;C_2,\]
where $C_1$ and $C_2$ are some symmetric matrices and  the components of $A$ and $B$ are given by $A_{ia,jb}=-\tfrac{1}{4}b(\gamma([X_i,X_j]),[X_a,X_b])$ and $B_{ij,ab}=- \tfrac{1}{2}b(\gamma([X_i,X_j]),[X_a,X_b])$, respectively. So $A$ is  a symmetric matrix, but one can think of it as being build from skew-symmetric blocks, arranged in a skew-symmetric way. 

It is then a simple exercise in linear algebra to show that such a matrix $A$ has a negative eigenvalue iff $A\neq0$.\\

The proof of proposition 3.3 
is now quite simple:

\begin{proof}[Proof of proposition 3.3] 
Let us first assume that $(M,g^M)$ has nonnegative curvature operator. Notice that $A=0$  if and only if $B=0$.
As $t \to 0$, the smallest eigenvalue of $\hat{R}^t$ converges to something which is less or equal to the smallest eigenvalue of $A$. Since $A$ has a negative eigenvalue if $A\neq 0$,  $(P,g^t)$ can only have almost nonnegative curvature operator as $t \to 0$ when $A=0$. But if $A=0$, then also $B=0$, so the smallest eigenvalue of $\hat{R}^t$ converges to 0 as $t \to 0$. It remains to note that $A=0$ if and only if $\mathrm{Im}(\Omega)$ is perpendicular to $[\mathfrak{g},\mathfrak{g}]$.

Let us now assume $M$ admits almost nonnegative curvature operator and that $\mathrm{Im}(\Omega)$ is perpendicular to $[\mathfrak{g},\mathfrak{g}]$, so $A=0$ for every metric on $M$. Let $\eps>0$ be given. We then choose a metric $g^M$  on $M$ such that $\hat{R}^M>-\tfrac{\eps}{2}$. In every point, the smallest eigenvalue of $\hat{R}^t$ then converges to the one of $\hat{R}^M$. By compactness, we can choose $t_0>0$ such that $\hat{R}^{t_0}>-\eps$ everywhere. 
\end{proof} 

We will now use Theorem 3.2 to prove Theorems 1 and 2
by looking at principal circle bundles over products of complex projective spaces.

\begin{proof}[Proof of Theorem 1 and Theorem 2]
Let now, for coprime nonzero integers $k$ and $l$, $P^{k,l}$ denote
the seven-dimensional total space of the principal circle bundle over $\CP^1 \times \CP^2$
with Euler class $k \alpha +  l \beta$, where $\alpha$ and $\beta$ are
generators of the second integral cohomology group of $\CP^1$ and $\CP^2$, respectively.
Then all $P^{k,l}$ are simply connected and admit ANCO by Theorem 3.2.

By an elementary calculation using the Gysin sequence one sees that
there are infinitely many distinct homotopy types
among the manifolds $P^{k,l}$, and by Remark~2.5
only finitely many of them can admit nonnegative curvature operator.
Crossing the others with simply connected spheres of appropriate codimension
then yields Theorem~1.

To obtain Theorem 2, we note that Kreck and Stolz (see [KS88]) showed
in particular that among the $P^{k,l}$ there exist homeomorphic manifolds
which are not diffeomorphic. (That actually none of these can admit 
nonnegative curvature operator is clear from Proposition~2.4.)
\end{proof}

\section{Further Remarks and Open Questions}

We conclude this article with a number of remarks and open questions.\\

The biggest question, of course, is whether it is possible to classify,
as in the case of nonnegative curvature operator, 
all manifolds which admit almost nonnegative curvature operator.\\

A positive answer to the following question
would, at least essentially, reduce the general classification to 
the simply connected case.

\noindent {\bf Question 4.1} \thinspace {\it Let $M$ be manifold
with almost nonnegative curvature operator. Is it true that a
finite cover $\tilde{M}$ of $M$ is the total space of a fibre
bundle
 $$
 F\to \tilde{M}\to N
 $$
over a nilmanifold $N$ with a simply connected fibre $F$ which
admits almost nonnegative curvature operator?}

Moreover, one may ask:

\noindent {\bf Question 4.2}  \thinspace {\it
Are there closed manifolds with almost nonnegative sectional curvature
which do not admit almost nonnegative curvature operator?}

\noindent {\bf Question 4.3} \thinspace {\it
Which (non-trivial) sphere bundles over spheres do admit almost nonnegative curvature operator?}

We would like to mention here that all
non-trivial $S^{4k+1}$ bundles over $S^2$, $k=1,2,\ldots$, admit ANCO.
They can be written as principal circle bundles 
over $S^2 \times \CP^{2k}$, so that Theorem 3.2 applies to them.

\noindent {\bf Question 4.4} \thinspace {\it
Are there any exotic spheres with almost nonnegative curvature operator?}

\noindent {\bf Question 4.5} \thinspace {\it
Are manifolds with almost nonnegative curvature operator rationally elliptic?}

\noindent {\bf Question 4.6} \thinspace {\it
Do manifolds with almost nonnegative curvature operator have nonnegative Euler characteristic?}

Notice that for manifolds with nonnegative curvature operator, 
the answer to the last two questions is positive.

\noindent {\bf Remark 4.7} \thinspace
Theorem 3.1 and Theorem 3.2 of course give rise to topologically quite
more complicated examples than the ones which we used to prove Theorem 1 and Theorem 2.
But none of them seem to solve the case of dimensions different from the
ones appearing there so far.

Nevertheless, we do suspect that the classes of simply connected
manifolds with nonnegative curvature operator and the respective one
of ANCO manifolds actually already differ from each other
in every dimension $n\ge 4$.

\

\small
\bibliographystyle{alpha}

\begin{thebibliography}{GVWZ06}



\bibitem[Be88]{Be88}
B\'erard, P. H. {\it From vanishing theorems to estimating theorems: the Bochner technique revisited.} Bull. Amer. Math. Soc. (N.S.) 19 (1988), no. 2, 371–-406.



\bibitem[BW08]{BW08}
B\"ohm, C.; Wilking, B. {\it Manifolds with positive curvature
operators are space forms.}
 Ann. of Math. (2) 167 (2008),  no. 3, 1079--1097.



\bibitem[CLN06]{CLN06} Chow, B.; Lu, P.; Ni, L.
{\it Hamilton's Ricci flow.}
Graduate Studies in Mathematics, 77. American Mathematical Society, Providence, RI; Science Press, New York, 2006.


\bibitem[FY92]{FY92}
Fukaya, K.; Yamaguchi, T.
\newblock {\it The fundamental groups of almost nonnegatively curved
manifolds.}
\newblock {Ann. of Math. (2)}, 136 (1992), no. 2, 253--333.



\bibitem[Hel01]{Hel01}
Helgason, S. {\it
 Differential geometry, Lie groups, and symmetric spaces.} Corrected reprint of the 1978 original. Graduate Studies in Mathematics, 34. American Mathematical Society, Providence, RI, 2001.


\bibitem[Her12]{Her12} Herrmann, M.
{\it Almost nonnegative curvature operator on certain principal bundles.}
In: Oberwolfach Reports, Report No. 01/2012,
Mini-Workshop: Manifolds with Lower Curvature Bounds
(organised by A. Dessai, W. Tuschmann, B. Wilking), 30--31, 2012.


\bibitem[KPT10]{KPT10}
Kapovitch, V.; Petrunin, A.; Tuschmann, W. {\em Nilpotency, almost
nonnegative curvature, and the gradient flow on Alexandrov
spaces.} Ann. Math. (2) 171 (2010), no. 1, 343--373.



\bibitem[KS88]{KS88}
Kreck, M.; Stolz, S. {\it
A diffeomorphism classification of 7-dimensional homogeneous Einstein manifolds with 
$SU(3) \times SU(2) \times U(1)$-symmetry.}
Ann. of Math. (2) 127 (1988), no. 2, 373-–388. 

\bibitem[Lo12]{Lo12}
Lott, J. {\it Collapsing with a lower bound on the curvature operator}. Preprint, 2012.


\bibitem[Se11]{Se11}
Sebastian, D.
\newblock {\it Konstruktionen von Mannigfaltigkeiten mit fast nichtnegativem
Kr\"ummungsoperator}. Dissertation, Karlsruher Institut f\"ur Technologie, 2011.


\bibitem[Tu12]{Tu12}
Tuschmann, W.
\newblock {\it Collapsing and Almost Nonnegative Curvature}.
In: C. B\"ar et al. (eds.), {\sl Global Differential Geometry}, pp. 93--106,
Springer Proceedings in Mathematics 17,
Springer, Berlin and Heidelberg, 2012.


\bibitem[Wil00]{Wil00} Wilking, B.
{\it 
On fundamental groups of manifolds of nonnegative curvature.}
Differential Geom. Appl. 13 (2000), no. 2, 129–-165.


\


\end{thebibliography}

\end{document}